\documentclass[11pt,reqno]{amsart}
\usepackage[T1]{fontenc}
\usepackage{lmodern}
\usepackage{amsmath,amssymb,amsthm,mathtools,mathrsfs}
\usepackage[margin=1.08in]{geometry}
\usepackage[hidelinks]{hyperref}
\usepackage{microtype}
\usepackage{enumitem}
\hypersetup{
  pdftitle={A Weak Lambda-Set That Is Not a Lambda-Set under p = c},
  pdfauthor={Alexander V. Osipov and Anton E. Lipin}
}
\newtheorem{theorem}{Theorem}[section]
\newtheorem{lemma}[theorem]{Lemma}
\newtheorem{proposition}[theorem]{Proposition}
\newtheorem{corollary}[theorem]{Corollary}
\theoremstyle{definition}

\newtheorem{problem}[theorem]{Problem}
\theoremstyle{remark}
\newtheorem{remark}[theorem]{Remark}
\newcommand{\Fin}{\mathrm{Fin}}
\newcommand{\FU}{\mathrm{FU}}
\newcommand{\supp}{\operatorname{supp}}

\newcommand{\bfrak}{\mathfrak{b}}
\newcommand{\pfrak}{\mathfrak{p}}
\newcommand{\tfrak}{\mathfrak{t}}
\newcommand{\cfrak}{\mathfrak{c}}
\newcommand{\om}{\omega}
\newcommand{\doi}[1]{\href{https://doi.org/#1}{doi:\nolinkurl{#1}}}

\title[Weak $\lambda$-sets under $\mathfrak p=\mathfrak c$]
{A Weak $\lambda$-Set That Is Not a $\lambda$-Set\\
under $\mathfrak p=\mathfrak c$}
\author[Alexander V. Osipov]{Alexander V. Osipov\textsuperscript{a,b}}
\author[Anton E. Lipin]{Anton E. Lipin\textsuperscript{a,b}}
\address[a]{N.\,N. Krasovskii Institute of Mathematics and Mechanics,\\
  Ural Branch of the Russian Academy of Sciences}
\address[b]{Ural Federal University}

\makeatletter
\renewcommand{\@setauthors}{%
  \begin{center}
    \def\and{\unskip,\nobreak\hspace{0.45em}\ignorespaces}%
    {\normalsize\authors\par}
    \medskip
    \small\itshape
    \def\author##1{}%
    \def\address##1##2{\textsuperscript{##1}\,##2\par\smallskip}%
    \addresses
  \end{center}%
  \global\let\addresses\@empty
}
\makeatother


\begin{document}
\begin{abstract}
Under $\mathfrak p=\mathfrak c$, Cantor space contains a weak
$\lambda$-set that is not a $\lambda$-set. This answers the
consistency question of de Oliveira Rodrigues, dos Santos Ronchim, and
Szeptycki. The hypothesis $\mathfrak p=\mathfrak c$ is compatible
with the failure of the continuum hypothesis. The construction
combines the finite-unions form of Hindman's theorem with Bell's
characterization of $\mathfrak p$ in terms of $\sigma$-centered
forcing. The increasing enumerations of the infinite members of the
constructed space form an unbounded family under eventual domination. This makes
the countable set of finite subsets of $\omega$ a witness to the
failure of the $\lambda$-property.
\end{abstract}
\maketitle

\section{Introduction and background}

A $\lambda$-set is a subspace of a Polish space in which every
countable subset is relatively $G_\delta$. Miller~\cite{Miller1984}
surveys this and other special sets of reals. For almost disjoint
families of branches, properties of the corresponding sets of reals
are closely related to separation properties of the associated
Isbell--Mr\'owka spaces~\cite{RodriguesRonchimSzeptycki2023}.

In 2023, de Oliveira Rodrigues, dos Santos Ronchim, and Szeptycki
\cite[Definition~4.1]{RodriguesRonchimSzeptycki2023} introduced the
notion of a \emph{weak $\lambda$-set}. In such a space any two disjoint
countable subsets can be separated by a relative $\Delta^0_2$ set. They
related this condition to strong $\aleph_0$-separation in almost
disjoint families of branches and proved that every weak $\lambda$-set
in a Polish space is perfectly meager
\cite[Theorem~4.4]{RodriguesRonchimSzeptycki2023}.
The first part of their Question~4.6 asks whether it is consistent that
there exists a weak $\lambda$-set that is not a $\lambda$-set.

Weak $\lambda$-spaces also arise in the study of spaces of Baire-one
functions with the topology of pointwise convergence.
Osipov~\cite{Osipov2026} relates the weak $\lambda$-property to the
Choquet property of these function spaces. Here we consider subsets
of Cantor space.

Our main result is the following conditional separation theorem.

\begin{theorem}\label{thm:main}
Assume $\pfrak=\cfrak$. There is a set $X\subseteq2^\om$ with the
following properties.
\begin{enumerate}[label=\textup{(\roman*)},leftmargin=2em]
\item For every pair of disjoint countable sets $A_0,A_1\subseteq X$,
      there is a partition $X=P_0\sqcup P_1$ into relative
      $\Delta^0_2$ sets with $A_i\subseteq P_i$ for $i=0,1$.
\item There is a countable subset of $X$ that is not $G_\delta$ in $X$.
\end{enumerate}
In particular, $X$ is a weak $\lambda$-set that is not a
$\lambda$-set.  The same conclusion holds for a subspace of the real
line.
\end{theorem}

The continuum hypothesis (CH) implies $\pfrak=\cfrak$. Bell's theorem
identifies $\pfrak=\cfrak$ with Martin's axiom for $\sigma$-centered
posets~\cite{Bell1981,Schilhan2019}. Thus the theorem also applies in
models of Martin's axiom in which $\cfrak>\aleph_1$.

ChatGPT (OpenAI) generated the proof of Theorem~\ref{thm:main}.
The proof constructs a decreasing chain of block sequences. At each stage, a
finite-unions lemma gives an $F_\sigma$ partition of a set containing
the finite subsets of $\om$ and all points chosen from that stage on.
At a limit stage, a chain of length less than $\pfrak$ has a common
lower bound. The earlier points form a relative $F_\sigma$ set of
cardinality less than $\bfrak$, so they form a $\lambda$-set even
when there are uncountably many of them.

\section{Notation and two cardinal invariants}

Let $\cfrak=2^{\aleph_0}$. Identify Cantor space $2^\om$ with
$\mathcal P(\om)$ via characteristic functions, and put
\[
\Fin=[\om]^{<\om},\qquad [\om]^\om=2^\om\setminus\Fin.
\]
For $a,b\subseteq\om$, write $a\subseteq^* b$ if $a\setminus b$
is finite. For $s\in2^{<\om}$, let $[s]$ denote the set of all
extensions of $s$ in $2^\om$.
A set is called $\Delta^0_2$ in a space if it is simultaneously
$F_\sigma$ and $G_\delta$ there. The complement of a relative
$\Delta^0_2$ set is again relative $\Delta^0_2$. Thus the separator
in \cite[Definition~4.1]{RodriguesRonchimSzeptycki2023} gives exactly
the partition in Theorem~\ref{thm:main}(i).

For $f,g\in\om^\om$, write $f\leq^*g$ if
$f(n)\leq g(n)$ for all but finitely many $n$, and write $f<^*g$
if $f(n)<g(n)$ eventually.  The bounding number $\bfrak$ is the
least cardinality of an unbounded subset of
$(\om^\om,\leq^*)$.  We use the standard inequalities
\[
\aleph_1\leq\pfrak\leq\bfrak\leq\cfrak.
\]
Recall that $\pfrak$ is the least size of a family
$\mathcal F\subseteq[\om]^\om$ with the strong finite intersection
property and no infinite pseudointersection. The former means that
every finite intersection is infinite. An infinite pseudointersection
of $\mathcal F$ is a set $a\in[\om]^\om$ such that
$a\subseteq^*F$ for all $F\in\mathcal F$.

We use the following consequence of Bell's characterization of
$\pfrak$~\cite{Bell1981,Schilhan2019}.

\begin{theorem}[Bell]\label{thm:bell}
If $\mathbb P$ is a nonempty $\sigma$-centered partially ordered set and
$\mathscr D$ is a collection of fewer than $\pfrak$ dense subsets of
$\mathbb P$, then there is a filter in $\mathbb P$ meeting every
member of $\mathscr D$.
\end{theorem}

The next lemma is the lower bound in
\cite[Proposition~4.7]{RodriguesRonchimSzeptycki2023}.

\begin{lemma}\label{lem:smalllambda}
Every subset $Y$ of $2^\om$ with $|Y|<\bfrak$ is a $\lambda$-set.
\end{lemma}
\begin{proof}
Let $A=\{a_i:i<\om\}$ be a countably infinite subset of $Y$,
with the $a_i$ distinct.  For each $y\in Y\setminus A$ define
$h_y\in\om^\om$ by choosing $h_y(i)$ so that
$y\notin[a_i\upharpoonright h_y(i)]$.  As there are fewer than
$\bfrak$ functions $h_y$, choose $g\in\om^\om$ with
$h_y\leq^*g$ for all $y\in Y\setminus A$.  Put
\[
U_n=\bigcup_{i<\om}
[a_i\upharpoonright\max\{g(i),n\}],\qquad n<\om.
\]
Every $U_n$ is open and contains $A$.  Fix $y\in Y\setminus A$.
For all sufficiently large $i$ we have $h_y(i)\leq g(i)$, and
for the remaining finitely many $i$ we may choose $n$ sufficiently
large that $h_y(i)\leq n$.  Then $y\notin U_n$.  Consequently
$A=Y\cap\bigcap_n U_n$.  Finite $A$ are $G_\delta$ in the metric
space $Y$ as well.
\end{proof}

\begin{lemma}[The finite-set test]\label{lem:finite-test}
Let $S\subseteq[\om]^\om$ and set $X=\Fin\cup S$.  For
$x\in[\om]^\om$ let $e_x\in\om^\om$ be its strictly increasing
enumeration.  Then
\[
\Fin\text{ is a relative }G_\delta\text{ in }X
\quad\Longleftrightarrow\quad
\{e_x:x\in S\}\text{ is }\leq^*\text{-bounded}.
\]
\end{lemma}
\begin{proof}
Suppose $\Fin=X\cap G$ with $G\subseteq2^\om$ a $G_\delta$.
Then $S\subseteq2^\om\setminus G=\bigcup_m K_m$, where the
$K_m$ are compact subsets of $[\om]^\om$.
For each fixed $n$ and $m$, compactness supplies an integer
$g_m(n)$ such that every $x\in K_m$ has at least $n+1$ elements
below $g_m(n)+1$.  Thus the enumerations from $K_m$ are
pointwise bounded by $g_m$. The function
$g(n)=\max_{m\le n}g_m(n)$ eventually dominates every $g_m$ and
hence every $e_x$ with $x\in S$.

Conversely, let $g\in\om^\om$ eventually dominate each $e_x$
for $x\in S$.  For each $m$ define
\[
K_m=\bigcap_{n\geq m}
\{z\in2^\om:|z\cap[0,g(n)]|\geq n+1\}.
\]
Each set in the intersection is clopen, so $K_m$ is compact.
Also $K_m\cap\Fin=\varnothing$, and
$S\subseteq\bigcup_m K_m$.  Therefore
$G=2^\om\setminus\bigcup_mK_m$ is a $G_\delta$ with
$X\cap G=\Fin$.
\end{proof}

\section{Block sequences and their \texorpdfstring{$F_\sigma$}{F-sigma} hulls}

A \emph{block sequence} is an infinite sequence
$B=(b_n)_{n<\om}$ of nonempty finite subsets of $\om$ satisfying
$\max b_n<\min b_{n+1}$.  Write
\[
\FU(B)=\left\{\bigcup_{n\in F}b_n:
\varnothing\ne F\in[\om]^{<\om}\right\},
\qquad \supp(B)=\bigcup_{n<\om}b_n.
\]
We call $C$ a \emph{condensation} of $B$, written
$C\preceq B$, if every block of $C$ belongs to $\FU(B)$.
Write $C\preceq^*B$ if all but finitely many blocks of $C$
belong to $\FU(B)$. Both relations are transitive. Indeed, blocks
are pairwise disjoint, so in a further condensation only finitely
many blocks can meet any of the finitely many exceptional blocks.
Define
\[
K(B)=\left\{\bigcup_{n\in I}b_n:I\subseteq\om\right\},
\quad
D(B)=\{s\mathbin\triangle y:s\in\Fin,\ y\in K(B)\},
\]
where $\triangle$ denotes symmetric difference.

\begin{lemma}\label{lem:hull}
Every block sequence $B$ has the following properties.
\begin{enumerate}[label=\textup{(\roman*)},leftmargin=2em]
\item $K(B)$ is homeomorphic to $2^\om$.
\item $D(B)$ is $F_\sigma$ in $2^\om$ and contains $\Fin$.
\item If $C\preceq^*B$, then $D(C)\subseteq D(B)$.
\item If $x\in D(B)$, then
$|x\setminus\supp(B)|<\aleph_0$.
\end{enumerate}
\end{lemma}
\begin{proof}
The map $I\mapsto\bigcup_{n\in I}b_n$ is a continuous
injection of the compact Cantor space into itself, and hence
is a homeomorphism onto its image $K(B)$.  Each translate
$s\triangle K(B)$ is compact, so their countable union
$D(B)$ is $F_\sigma$.  Taking $y=\varnothing$ shows
$\Fin\subseteq D(B)$.  If $C\preceq^*B$, any union of
blocks of $C$ differs by a finite set from a union of blocks
of $B$, proving (iii).  Finally, all members of $K(B)$ are
subsets of $\supp(B)$, proving (iv).
\end{proof}

\section{The finite-unions separation lemma}

By the finite-unions form of Hindman's theorem~\cite{Hindman1974},
every finite coloring of $\FU(C)$ admits a condensation $B\preceq C$
such that $\FU(B)$ is monochromatic. We use this theorem to extend a
coloring of $\Fin$ to a partition of $D(B)$ into two $F_\sigma$ sets.

\begin{lemma}[Finite-unions separation]\label{lem:sep}
Let $C$ be a block sequence and let
$c:\Fin\to\{0,1\}$ be a coloring.  There exist a condensation
$B\preceq C$ and disjoint $F_\sigma$ subsets
$H_0,H_1\subseteq2^\om$ such that
\[
D(B)=H_0\sqcup H_1,
\qquad H_i\cap\Fin=c^{-1}(i)\quad(i=0,1).
\]
\end{lemma}
\begin{proof}
We construct $B=(b_k)_{k<\om}\preceq C$ with the property that,
for each $k$ and every $t\subseteq[0,\max b_k]$, there is a
color $\varepsilon_{k,t}\in\{0,1\}$ such that
\begin{equation}\label{eq:hom}
c(t\cup u)=\varepsilon_{k,t}
\quad\text{for all }u\in\FU((b_n)_{n>k}).
\end{equation}
Take $b_0$ to be the first block of $C$ and let $T_0$ be its tail.
Assume that $b_k$ and a block tail $T_k$ to its right have
been selected.  Color $\FU(T_k)$ by the finite vector
\[
u\longmapsto \bigl(c(t\cup u)\bigr)_
{t\subseteq[0,\max b_k]}.
\]
Hindman's theorem gives $S_k\preceq T_k$ on whose finite
unions this vector is constant.  Let $b_{k+1}$ be the first
block of $S_k$ and let $T_{k+1}$ be its remaining tail.
Every $b_n$ with $n>k+1$ is a finite union of blocks of $T_{k+1}$.
Together with $b_{k+1}$, these blocks therefore form a condensation of
$S_k$. All their nonempty finite unions lie in $\FU(S_k)$, which gives
\eqref{eq:hom}. The construction also gives $B\preceq C$.

For $k<j$ and $t\subseteq[0,\max b_k]$, put
\[
L_{k,t,j}=
\left\{t\cup b_j\cup\bigcup_{n\in I}b_n:
I\subseteq\{j+1,j+2,\ldots\}\right\}.
\]
Each $L_{k,t,j}$ is a continuous image of the Cantor space of choices
for $I$, and is therefore compact. Its elements differ from unions of
blocks of $B$ only on the finite set $t$, so $L_{k,t,j}\subseteq D(B)$.
Every finite point of $L_{k,t,j}$ has the form $t\cup u$ for
some $u\in\FU((b_n)_{n>k})$, and hence has the fixed color
$\varepsilon_{k,t}$. Add every singleton $\{s\}$ with $s\in\Fin$
to this family and assign it the color $c(s)$.

These compact sets form a countable cover of $D(B)$. Indeed,
if $x\in D(B)\setminus\Fin$, write
$x=s\triangle\bigcup_{n\in I}b_n$ where $s$ is finite and
$I$ is infinite.  Choose $k$ with
$s\subseteq[0,\max b_k]$, put
$t=x\cap[0,\max b_k]$, and let $j$ be the least member
of $I$ larger than $k$.  Then $x\in L_{k,t,j}$.

Enumerate the compact cover as $(L_m)_{m<\om}$, retaining
for each $L_m$ a color $\delta_m$ consistent with $c$ on
$L_m\cap\Fin$.  Put
\[
E_m=L_m\setminus\bigcup_{\ell<m}L_\ell,
\qquad
H_i=\bigcup_{\delta_m=i}E_m.
\]
Each $E_m$ is the intersection of a closed set with an open set.
Since open subsets of Cantor space are $F_\sigma$, so is each $E_m$.
The $E_m$ are disjoint and cover $D(B)$. If $s\in E_m\cap\Fin$,
then $s\in L_m$ and $\delta_m=c(s)$. This proves the assertion
about $H_0$ and $H_1$.
\end{proof}

\section{Common lower bounds below \texorpdfstring{$\mathfrak p$}{p}}

\begin{lemma}[The $\pfrak$-diagonalization lemma]\label{lem:pdiag}
Suppose $\kappa$ is an ordinal with $|\kappa|<\pfrak$ and
$(B_\xi)_{\xi<\kappa}$ is a family of block sequences satisfying
\[
B_\eta\preceq^*B_\xi\qquad(\xi<\eta<\kappa).
\]
Then there exists a block sequence $C$ such that
$C\preceq^*B_\xi$ for every $\xi<\kappa$.
\end{lemma}
\begin{proof}
For $\kappa=0$, take $C=(\{n\})_{n<\om}$. Otherwise let
$\mathbb P$ consist of pairs $p=(s,F)$, where $s$ is a finite
block sequence and $F\in[\kappa]^{<\om}$.  We say that
$(t,G)\leq(s,F)$ if $s$ is an initial segment of $t$,
$F\subseteq G$, and every block of $t$ beyond $s$ belongs to
\[
\bigcap_{\xi\in F}\FU(B_\xi).
\]
Fix a stem $s$. Any finite family $(s,F_1),\ldots,(s,F_r)$ has
the common extension $(s,F_1\cup\cdots\cup F_r)$. Since there
are only countably many finite stems, $\mathbb P$ is $\sigma$-centered.

For each $\xi<\kappa$, the set
$E_\xi=\{(s,F):\xi\in F\}$ is dense.  For each $n<\om$,
$D_n=\{(s,F):|s|\geq n\}$ is also dense.  To see this,
fix $(s,F)$ with $F\ne\varnothing$ and let
$\eta=\max F$. For each $\xi\in F$, all but finitely many blocks of
$B_\eta$ belong to $\FU(B_\xi)$. Since $F$ is finite, every
sufficiently late block belongs to all these sets. We can append
as many of these blocks as needed beyond the last block of $s$.
For $F=\varnothing$ we may append any blocks to the right of $s$.

There are at most $|\kappa|+\aleph_0<\pfrak$ dense sets in this
family. By Theorem~\ref{thm:bell}, choose a filter meeting all of
them. Its stems are linearly ordered by extension, and meeting
every $D_n$ makes their union an infinite block sequence $C$.
Fix $(s,F)$ in the filter with $\xi\in F$. For any block of $C$
beyond $s$, choose a condition in the filter whose stem includes
that block. A common extension of the two conditions shows that the
block belongs to $\FU(B_\xi)$. Thus $C\preceq^*B_\xi$ for every
$\xi<\kappa$.
\end{proof}

\begin{lemma}[Sparse refinement]\label{lem:sparse}
Let $B=(b_n)_{n<\om}$ be a block sequence and let
$g\in\om^\om$.  There exists a block subsequence
$C=(b_{k_n})_{n<\om}$ such that, writing
$x=\supp(C)$, one has $e_x\not\leq^*g$.
\end{lemma}
\begin{proof}
Choose $k_0<k_1<\dots$ recursively.  If
$k_0,\ldots,k_{n-1}$ have been chosen, let
$m_n=\sum_{j<n}|b_{k_j}|$ and pick $k_n$ after the previous
index so that $\min b_{k_n}>g(m_n)$.  Since the blocks are
separated and their minima tend to infinity, this is
possible.  The $m_n$-th element of $x$ (numbered starting
at zero) equals $\min b_{k_n}$, so
$e_x(m_n)>g(m_n)$ for all $n$.  The indices $m_n$ tend to
infinity, as required.
\end{proof}

\section{Increasing enumerations and initial segments}

The supports will decrease modulo finite sets. Their enumerations
allow us to separate the earlier points from the rest of the space.

\begin{lemma}\label{lem:strict}
Suppose $A,B\in[\om]^\om$, $A\subseteq^*B$ and
$|B\setminus A|=\aleph_0$.  Then $e_B<^*e_A$.
\end{lemma}
\begin{proof}
The set $A\setminus B$ is finite while $B\setminus A$ is
infinite.  Consequently, the difference
$|B\cap[0,m]|-|A\cap[0,m]|$ tends to infinity with $m$.
For all sufficiently large $n$, taking $m=e_A(n)$ gives at least
$n+2$ elements of $B$ in $[0,e_A(n)]$. Hence $e_B(n)<e_A(n)$.
\end{proof}

\begin{lemma}[Initial-segment lemma]\label{lem:initial}
Let $\kappa$ be any ordinal and
$S=\{x_\alpha:\alpha<\kappa\}\subseteq[\om]^\om$
be a family with
\[
e_{x_\alpha}<^*e_{x_\beta}\qquad(\alpha<\beta<\kappa).
\]
Set $X=\Fin\cup S$.  For each $\alpha<\kappa$, the
initial segment
$S_{<\alpha}=\{x_\beta:\beta<\alpha\}$ is
relatively $F_\sigma$ in $X$, as is its complement.
In particular, both are relative $\Delta^0_2$ sets.
\end{lemma}
\begin{proof}
Write $f_\alpha=e_{x_\alpha}$ and define
\[
F_\alpha=\bigcup_{m<\om}\bigcap_{n\geq m}
\left\{z\subseteq\om:
|z\cap[0,f_\alpha(n))|\geq n+1\right\}.
\]
Each set in braces is clopen, so $F_\alpha$ is $F_\sigma$ in
$2^\om$. No finite set belongs to $F_\alpha$. For infinite $z$,
membership means that $e_z(n)<f_\alpha(n)$ eventually. The strict
increase of the family gives
\[
X\cap F_\alpha=S_{<\alpha}.
\]
For the complementary set use
\[
V_\alpha=\bigcup_{m<\om}\bigcap_{n\geq m}
\left\{z\subseteq\om:
|z\cap[0,f_\alpha(n))|\leq n\right\}.
\]
Again this is $F_\sigma$.  Every finite set belongs to
$V_\alpha$, and for infinite $z$ membership means
$f_\alpha\leq^*e_z$.  Consequently,
\[
X\cap V_\alpha=\Fin\cup\{x_\beta:\alpha\leq\beta<\kappa\}
=X\setminus S_{<\alpha}.
\]
Both complementary pieces are $F_\sigma$ in $X$ and
hence $\Delta^0_2$ there.
\end{proof}

\section{Construction and proof of the main theorem}

\begin{proof}[Proof of Theorem~\ref{thm:main}]
Assume $\pfrak=\cfrak$, so that also
$\bfrak=\cfrak$. Enumerate all pairs
$(g,c)\in\om^\om\times\{0,1\}^{\Fin}$ as
$(g_\alpha,c_\alpha)_{\alpha<\cfrak}$. Each function $g$ occurs,
and each coloring $c$ occurs at $\cfrak$ many indices, hence
cofinally often.

By transfinite recursion on $\alpha<\cfrak$ we construct
block sequences $B_\alpha$, points
$x_\alpha=\supp(B_\alpha)$, and disjoint $F_\sigma$
sets $H_{\alpha,0},H_{\alpha,1}\subseteq2^\om$ satisfying
\begin{align}
B_\alpha&\preceq^*B_\beta
&& (\beta<\alpha), \label{eq:rec-chain}\\
|x_\beta\setminus x_\alpha|&=\aleph_0
&& (\beta<\alpha), \label{eq:rec-diff}\\
e_{x_\alpha}&\not\leq^*g_\alpha,
&& \label{eq:rec-unbounded}\\
D(B_\alpha)&=H_{\alpha,0}\sqcup H_{\alpha,1},
\qquad H_{\alpha,i}\cap\Fin=c_\alpha^{-1}(i)\quad(i=0,1).
&& \label{eq:rec-colors}
\end{align}

At stage $\alpha>0$, we have $|\alpha|<\cfrak=\pfrak$.
Lemma~\ref{lem:pdiag} gives a block sequence $C=(d_n)$ with
$C\preceq^*B_\beta$ for all $\beta<\alpha$. At stage $0$,
take $d_n=\{n\}$. Let
$C^{\mathrm{ev}}=(d_{2n})_{n<\om}$.
For each $\beta<\alpha$, all but finitely many blocks of
$C$ lie inside $x_\beta$. Since the blocks of $C$ are pairwise
disjoint, infinitely many odd-indexed blocks lie entirely in
$x_\beta\setminus\supp(C^{\mathrm{ev}})$. Hence
\begin{equation}\label{eq:discard}
|x_\beta\setminus\supp(C^{\mathrm{ev}})|=\aleph_0.
\end{equation}
Apply Lemma~\ref{lem:sep} to $C^{\mathrm{ev}}$ and the
coloring $c_\alpha$.  We obtain a condensation $B'$
and a partition
$D(B')=H'_0\sqcup H'_1$ into $F_\sigma$ sets with the
specified colors on $\Fin$.

Next use Lemma~\ref{lem:sparse} to obtain a block
subsequence $B_\alpha\preceq B'$ whose support
$x_\alpha$ satisfies
$e_{x_\alpha}\not\leq^*g_\alpha$.  Define
$H_{\alpha,i}=H'_i\cap D(B_\alpha)$.
By Lemma~\ref{lem:hull}, these sets are $F_\sigma$ and
form a partition of $D(B_\alpha)$, with the correct colors
on $\Fin$.  Transitivity of block refinement proves
\eqref{eq:rec-chain}.  Since
$x_\alpha\subseteq\supp(C^{\mathrm{ev}})$, \eqref{eq:discard}
proves \eqref{eq:rec-diff}.  The remaining requirements
follow by construction.

Set
\[
S=\{x_\alpha:\alpha<\cfrak\},
\qquad X=\Fin\cup S.
\]
For $\beta<\alpha$, \eqref{eq:rec-chain} and \eqref{eq:rec-diff}
give $x_\alpha\subseteq^*x_\beta$ and
$|x_\beta\setminus x_\alpha|=\aleph_0$.  Thus the
$x_\alpha$ are distinct, and Lemma~\ref{lem:strict} yields
\begin{equation}\label{eq:scale}
e_{x_\beta}<^*e_{x_\alpha}
\qquad(\beta<\alpha<\cfrak).
\end{equation}
For every $g\in\om^\om$, some $g_\alpha$ equals $g$, and
\eqref{eq:rec-unbounded} gives $e_{x_\alpha}\not\leq^*g$.
Thus $\{e_{x_\alpha}:\alpha<\cfrak\}$ is unbounded.
By Lemma~\ref{lem:finite-test}, $\Fin$ is \emph{not} a
relative $G_\delta$ subset of $X$.  Therefore $X$ is not
a $\lambda$-set.

It remains to verify weak $\lambda$-separation.  First,
Lemma~\ref{lem:hull} and the recursion give the exact tail
identity
\begin{equation}\label{eq:tail}
X\cap D(B_\alpha)
=\Fin\cup\{x_\beta:\alpha\leq\beta<\cfrak\}.
\end{equation}
Indeed, for $\beta\geq\alpha$ we have
$B_\beta\preceq^*B_\alpha$, hence
$x_\beta\in D(B_\beta)\subseteq D(B_\alpha)$.
For $\beta<\alpha$, condition
\eqref{eq:rec-diff} and Lemma~\ref{lem:hull}(iv) imply
$x_\beta\notin D(B_\alpha)$.

Let $A_0,A_1\subseteq X$ be disjoint countable sets.
Choose a coloring $c:\Fin\to2$ satisfying
$c(s)=i$ for $s\in A_i\cap\Fin$. By K\"onig's theorem,
$\operatorname{cf}(\cfrak)>\aleph_0$, so the countably many indices
of points in $(A_0\cup A_1)\cap S$ are bounded in $\cfrak$.
Choose $\alpha$ above all these indices such that
$c_\alpha=c$.

Put $E_\alpha=\{x_\beta:\beta<\alpha\}$.
By Lemma~\ref{lem:initial} and \eqref{eq:scale}, this is
an $F_\sigma$ subset of $X$.  Moreover,
$|E_\alpha|<\cfrak=\bfrak$, so
Lemma~\ref{lem:smalllambda} says that $E_\alpha$ is a
$\lambda$-set.  Write
\[
E_{\alpha,0}=A_0\cap E_\alpha,\qquad
E_{\alpha,1}=E_\alpha\setminus E_{\alpha,0}.
\]
The countable set $E_{\alpha,0}$ is $F_\sigma$ in $E_\alpha$,
and the $\lambda$-property makes it $G_\delta$ there. Hence
$E_{\alpha,1}$ is also $F_\sigma$ in $E_\alpha$. Since $E_\alpha$
is $F_\sigma$ in $X$, both pieces are $F_\sigma$ in $X$.

Now define, for $i=0,1$,
\[
P_i=(X\cap H_{\alpha,i})\cup E_{\alpha,i}.
\]
By \eqref{eq:tail} these two sets are disjoint and cover
$X$. They are complementary $F_\sigma$ sets in $X$, so both are
also $G_\delta$ there. For $i=0,1$, the coloring puts
$A_i\cap\Fin$ in $H_{\alpha,i}$, and the choice of $\alpha$ puts
$A_i\cap S$ in $E_{\alpha,i}$. The partition therefore separates
$A_0$ and $A_1$, as required.

Finally, the standard embedding of $2^\om$ as the
middle-thirds Cantor set transfers the example to
$\mathbb R$.
\end{proof}

\begin{corollary}\label{cor:MA}
$\mathrm{MA}(\sigma\text{-centered})$ implies that the real line
contains a weak $\lambda$-set that is not a $\lambda$-set.
In particular, the existence of such a set is consistent with
$\neg\mathrm{CH}$ relative to ZFC.
\end{corollary}
\begin{proof}
Bell's theorem gives $\pfrak=\cfrak$, so Theorem~\ref{thm:main}
applies. For the consistency assertion, use the relative consistency
of $\mathrm{MA}+\neg\mathrm{CH}$ with ZFC~\cite{SolovayTennenbaum1971}.
\end{proof}

\begin{remark}\label{rem:scale}
The sequence $(e_{x_\alpha})_{\alpha<\cfrak}$ is $<^*$-increasing,
unbounded, and of length $\bfrak=\cfrak$, so it is a $\bfrak$-scale.
Moreover, $S$ itself is a $\lambda$-set. Given a countable
$T\subseteq S$, choose $\alpha$ above all of its indices. The initial segment
$S_{<\alpha}$ is $G_\delta$ in $S$ by
Lemma~\ref{lem:initial}. Its cardinality is less than
$\bfrak$, so it is a $\lambda$-set by
Lemma~\ref{lem:smalllambda}.  Thus $T$ is $G_\delta$
in $S_{<\alpha}$ and therefore in $S$.
\end{remark}

\section{The role of the hypothesis \texorpdfstring{$\mathfrak p=\mathfrak c$}{p = c}}

The proof uses $\pfrak=\cfrak$ to obtain lower bounds at every
stage of a recursion of length $\cfrak$. It also uses
$\bfrak=\cfrak$ to make each proper initial segment a $\lambda$-set.
The bound in Lemma~\ref{lem:pdiag} is sharp for arbitrary decreasing
chains. To see this, recall that the tower number $\tfrak$ is the
least cardinal length of a $\subseteq^*$-decreasing sequence of
infinite subsets of $\om$ with no infinite pseudointersection.
Malliaris and Shelah proved that
$\pfrak=\tfrak$ in ZFC~\cite{MalliarisShelah2016}.

\begin{proposition}\label{prop:sharp}
In ZFC there is a $\preceq^*$-decreasing chain of $\pfrak$
block sequences with no common $\preceq^*$-lower bound.
\end{proposition}
\begin{proof}
By $\pfrak=\tfrak$, choose a tower
$(T_\alpha)_{\alpha<\pfrak}$ with
$T_\beta\subseteq^*T_\alpha$ when $\alpha<\beta$
and with no infinite pseudointersection.  Enumerate each
$T_\alpha$ increasingly as $(t_{\alpha,n})_{n<\om}$
and form the singleton block sequence
$B_\alpha=(\{t_{\alpha,n}\})_{n<\om}$.
Then $B_\beta\preceq^*B_\alpha$ for
$\alpha<\beta$.  If $C$ were a common lower bound, its
infinite support would satisfy
$\supp(C)\subseteq^*T_\alpha$ for all $\alpha$, contrary
to the choice of the tower.
\end{proof}

There are $\cfrak$ colorings of $\Fin$. When $\pfrak<\cfrak$, a
recursion through all of them reaches stages where
Lemma~\ref{lem:pdiag} no longer applies. Proposition~\ref{prop:sharp}
shows that a decreasing chain can already fail to have a lower bound
at length $\pfrak$. This leaves open the possibility of choosing a
chain with additional structure or using a different construction.

\begin{problem}[The ZFC problem]\label{prob:zfc}
Does ZFC prove that some $X\subseteq 2^\om$ is a weak $\lambda$-set
but not a $\lambda$-set?
\end{problem}

\section*{Use of generative AI}
ChatGPT (OpenAI) generated the proof of the main theorem, including
its transfinite construction. The authors also used ChatGPT during
revision to review the mathematical arguments and edit the English.
The authors take full responsibility for the mathematical content
and the final text.

\end{document}